\documentclass{article}
\usepackage{amsmath,amssymb,enumerate}

\newtheorem{prop}{Proposition}
\newtheorem{coro}[prop]{Corollaire}
\newtheorem{lem}[prop]{Lemme}

\newtheorem{defi}[prop]{D\'efinition}
\newtheorem{theo}[prop]{Th\'eor\`eme}

\newcommand{\R}{\mathbb R}
\newcommand{\Q}{\mathbb Q}
\newcommand{\C}{\mathbb C}
\newcommand{\N}{\mathbb N}

\newcommand{\Z}{\mathbb Z}
\newcommand{\X}{\mathcal X}

\newcommand{\wt}{\widetilde}

\newenvironment{prv}{Preuve:}{$\Box$}

\begin{document}

\title{Sur l'uniformisation des orbifolds K\"ahl\'eriens compacts}
\author{ Philippe Eyssidieux}
\date{3 F\'evrier 2013}

\maketitle

Cet article g\'en\'eralise aux orbifolds k\"ahl\'eriens les r\'esultats g\'en\'eraux 
sur l'uniformisation des vari\'et\'es k\"ahl\'eriennes compactes pr\'esent\'es par l'auteur
au premier congr\`es joint des Soci\'et\'es Math\'ematiques Vietnamiennes et Fran\c{c}aises
organis\'e \`a Hu\^e du 20 au 24 Ao\^ut 2012.

\section{Conjecture de Shafarevich champ\^etre}
\subsection{Uniformisation champ\^etre}
Soit
$\X $ un orbifold K\"ahl\'erien connexe (s\'epar\'e). 
Soit $\gamma:\X  \to X$ le morphisme 
canonique vers son espace de modules qui est un espace analytique normal \`a singularit\'es quotient, propre si $\X$ l'est. 
Par exemple, $\X $ peut \^etre l'analytifi\'e d'un champ alg\'ebrique de Deligne Mumford s\'epar\'e et lisse sur $\C$ auquel cas $X$ est l'analytifi\'e d'un espace alg\'ebrique, c'est \`a dire un espace de Moishezon si $\X$ est propre. 
On fait de plus l'hypoth\`ese  que $\X$ n'a pas de groupe d'inertie au
point g\'en\'erique, de sorte que $\gamma$ est une \'equivalence au dessus d'un ouvert dense de $X$. 

Par un abus de langage, on identifiera 
$\X$ et $\X^{Betti}$ le champ topologique correspondant qui est de Deligne-Mumford \`a groupes
 d'inertie finis\footnote{Nous utilisons la r\'ef\'erence commode \cite{No1,No2} pour les champs topologiques.}.  Soit $ x\in\X(\C)$ 
un point complexe et $\pi_1(\X, x)$ le groupe fondamental de $\X$ . 
Pour $x,y\in X(\C)$,  toute application continue $\gamma:([0,1],\{ 0, 1\}) \to (\X,  x,  y)$ d\'efinit un
isomorphisme $\gamma_*: \pi_1(\X, x)\to\pi_1(\X, y)$. De plus si $I_y:=\pi^{loc}_1(\X,y)$ est le groupe d'inertie en $\bar y$ (qui est fini)
on a un morphisme de groupe naturel $ l_y:\pi^{loc}_1(\X,y)\to \pi_1(\X,y)$ puis
 $\gamma_*^{-1} \circ l_y:\pi^{loc}_1(\X,y)\to \pi_1(\X,x)$. Ces morphimes de groupes seront appel\'es les morphismes
locaux d'inertie \cite{No1}. 
Le groupe fondamental \'etale de $\X$ est par d\'efinition le compl\'et\'e profini $\pi_1^{et}(\X,x)$ de $\pi_1(\X,x)$
et le compos\'e $ l^{et}_y:\pi^{loc}_1(\X,y)\to \pi_1^{et}(\X,y)$ de 
$ l_y:\pi^{loc}_1(\X,y)\to \pi_1(\X,y)$ avec la surjection canonique $ l: \pi_1(\X,y)\to \pi_1^{et}(\X,y)$
sera appel\'e le morphisme \'etale-local d'inertie en $y$. 

\begin{defi}
 On dit que $\X$ est  d\'eveloppable (ou uniformisable  en topologie transcendante)
si tous les morphismes locaux d'inertie sont injectifs, resp.
uniformisable (ou uniformisable en topologie \'etale) si
 tous les morphismes \'etale-locaux d'inertie sont injectifs.
\end{defi}

 Le champ propre $\X$ est uniformisable si et seulement si il existe un espace analytique normal connexe propre et lisse
 $U$,  un groupe $G$ fini et 
une action de $G$ sur $U$ telle que  $ \X = [U/ G]$. Si l'espace des modules de $\X$ est 
projectif-alg\'ebrique, resp. K\"ahl\'erien, resp  quasi-projectif,  on voit ais\'ement que 
$U$ est \'egalement projectif-alg\'ebrique, resp.  K\"ahl\'erien, resp. quasi-projectif.

Si $\X$ est d\'eveloppable, on dispose d'une vari\'et\'e
 complexe connexe et simplement connexe $\wt{X^{u}}$, d'une action proprement discontinue de $\pi_1(\X,x)$ sur  $\wt{X}$ et d'un isomorphisme $\X \to [\pi_1(\X,x) \backslash \wt{X^{u}} ].$ 
Si $\X$ est uniformisable, \'ecrivant $\X = [U/ G]$ comme ci-dessus, 
$\wt{X^{u}}$ s'identifie au rev\^etement universel de $U$. 

Les deux notions uniformisable et uniformisable en topologie transcendante
 ne sont pas \'equivalentes a priori mais nous ne savons pas s'il est possible de construire un champ de Deligne-Mumford
propre et lisse 
uniformisable en topologie transcendante sans \^etre uniformisable. 

\subsection{Probl\`eme de Serre champ\^etre}

Le probl\`eme de caract\'eriser la classe des groupes k\"ahl\'eriens 
(groupes de pr\'esenta\- tion finie pouvant apparaitre comme un groupe fondamental de la forme $\pi_1(X,x)$
avec $X$ vari\'et\'e k\"ahl\'erienne compacte) est connu sous le nom de probl\`eme de Serre. 
Serre a montr\'e que tout groupe fini est k\"ahl\'erien mais il y a des obstructions pour les groupes infinis. Par exemple
l'ab\'elianis\'e d'un groupe k\"ahl\'erien est de rang impair par la th\'eorie de Hodge. 
Le probl\`eme de Serre s'\'etend naturellement \`a caract\'eriser les groupes k\"ahl\'eriens
orbifolds (groupes de pr\'esentation finie pouvant apparaitre comme un groupe fondamental de la forme $\pi_1(\X,x)$
avec $\X$ orbifold k\"ahl\'erien compact). Je ne connais pas d'obstruction \`a ce que les restrictions connues
 portant sur les groupes k\"ahl\'eriens (voir \cite{ABCKT} pour un survey) s'\'etendent 
aux groupes k\"ahl\'eriens orbifold. 
Une instance particuli\`ere du probl\`eme de Serre est la conjecture de Toledo 
qui pr\'edit qu'un groupe k\"ahl\'erien  a un second nombre de Betti positif. 
Il est donc tentant d'\'etendre la conjecture de Toledo
aux groupes k\"ahl\'eriens orbifolds.

La classe des groupes de pr\'esentation finie pouvant apparaitre comme un groupe fondamental de la forme $\pi_1(\X,x)$
avec $\X$ l'analytifi\'e d'un champ alg\'ebrique sur $\C$ s\'epar\'e, propre, lisse et de Deligne-Mumford
(groupes DM-projectifs) est a priori plus vaste que la classe des groupes fondamentaux des vari\'et\'es projectives lisses
ou d'espaces alg\'ebriques propres et lisses (groupes projectif-alg\'ebriques).

Il n'est pas connu si la classe des groupes k\"ahl\'eriens orbifolds  est plus vaste que celle 
des groupes DM-projectifs, ni si celle des groupes k\"ahl\'erien est plus vaste que 
celle des groupes projectifs-alg\'ebriques. Cependant:

\begin{lem}
 Si $\X$ est uniformisable (resp. un champ alg\'ebrique de Deligne-Mumford uniformisable), son groupe fondamental est k\"ahl\'erien (resp. projectif-alg\'ebrique). 
\end{lem}

\begin{prv}
 Soit $U\to \X$ un rev\^etement \'etale propre comme plus haut. Soit $G=\pi_1(\X)/ \pi_1(U)$ le quotient fini correspondant. 
Soit $S$ une surface projective lisse telle que $\pi_1(S)=G$.
Soit $S'$ son rev\^etement universel. Le groupe $\pi_1(\X,x)$ agit sur $S'$  via son quotient $G$ et
son action diagonale sur $\wt{X^{u}} \times S'$ est propre et sans  point fixe. Par suite, posant 
$X'=\pi_1(\X,x) \backslash \wt{X^{u}} \times S'$,  $\pi_1(X')=\pi_1(\X)$. 
\end{prv}

\subsection{Convexit\'e holomorphe}

\begin{defi}
Un espace complexe normal $S$ est  holomorphiquement convexe
 s'il existe une application holomorphe propre \`a fibres connexes $\pi : S\to T$  
telle que $T$  est un espace de Stein normal. On appelle $T=Red(S)$
 la r\'eduction de Cartan-Remmert  de $S$.
\end{defi}

La conjecture de Shafarevich sur l'uniformisation (voir \cite{pasy} pour une introduction)
stipule que le rev\^etement universel d'une vari\'et\'e projective lisse est holomorphiquement 
convexe.

La conjecture de Shafarevich  est invariante bim\'eromorphiquement. Par suite, la convexit\'e holomorphe du rev\^etement 
universel d'un espace alg\'ebrique propre et lisse
r\'esulte de celle du rev\^etement 
universel d' un mod\`ele birationnel  projectif. 
De m\^eme, elle  
 implique aussi que si $\X$ est un champ de Deligne-Mumford lisse propre et 
uniformisable alors $\wt{\X^{u}}$ est holomorphiquement convexe. 
En revanche, le cas K\"ahl\'erien ne r\'esulte pas du cas projectif.

La conjecture de Shafarevich est ouverte. Une  classe d'exemples o\`u elle n'est pas \`a ce jour d\'ecid\'ee est celle des surfaces projectives fibr\'ees en courbes de genre 2. 
Plus g\'en\'eralement, on \'etudie la question suivante
\footnote{\`A partir de ce point on suppose $\X$ compact. }:
  sous quelles conditions portant sur $(\X, H)$ 
avec  $H\subset \pi_1(\X,x)$  un sous groupe distingu\'e
peut on affirmer que $H\backslash \wt{\X^{u}}$ est holomorphiquement convexe? Si cel\`a a lieu on dit que
$(\X, H)$ v\'erifie (HC)
.

\begin{lem} \label{shafmor}
 Soit $(\X, H)$ v\'erifiant (HC). Le groupe $\Gamma_H=\pi_1(\X,x)/H$ op\'ere proprement discontinuement sur 
$Red(H\backslash \wt{\X^{u}})$, l'application naturelle $H\backslash \wt{\X^{u}}\to Red(H\backslash \wt{\X^{u}})$
est $\Gamma_H$-\'equivariante et
d\'efinit une application holomorphe de champs analytiques 
$\mathfrak{s}^H: \X^{an} \to [ \Gamma_H \backslash Red(H\backslash \wt{\X^{u}}) ]$. 
\end{lem}

\begin{coro} \label{bezout}
 Soit $(\X, H)$ v\'erifiant (HC) et tel que toute application holomorphe 
d\'efinie sur $X$ est finie ou constante. Alors
ou bien $\Gamma_H$ est fini ou bien $H \backslash\wt{\X^{u}}$ est une vari\'et\'e de Stein. 
\end{coro}

\begin{prv}
 Passant \`a l'espace grossier $\mathfrak{s}^H$ induit une application holomorphe
$ s: X \to s(X)= \Gamma \backslash Red(H\backslash \wt{\X^{u}})$. 
L'hypoth\`ese implique que ou bien $s$ est constante ou bien $s$ a toutes ses fibres de dimension $0$.
En particulier $H\backslash\wt{\X^{u}} \to Red(H\backslash \wt{\X^{u}})$ est ou bien constant (auquel cas $H\backslash\wt{\X^{u}}$ est
compact et $\Gamma$ fini) ou bien \`a fibre finies et connexes c'est \`a dire un isomorphisme d'espaces normaux. 
\end{prv}

L'hypoth\`ese est v\'erifi\'ee par exemple si $X$ est un espace projectif. 

\begin{coro}
 Soit $\X$ d\'eveloppable tel que  $h^{1,1}(X)=1$. Supposons que $\pi_1(\X)$ n'a pas de quotient infini violant la propri\'et\'e de Kazhdan.  
 Alors
 $\wt{\X^{u}}$ est une vari\'et\'e de Stein. 
\end{coro}

\begin{prv} L'hypoth\`ese est qu'il existe $H\subset \pi_1(\X)$ distingu\'e tel 
 que  $\Gamma=H\backslash  \pi_1(\X)$ 
 a un module hilbertien $L^2$ tel que $H^1(\Gamma, L^2)\not =0$.  Mok a construit sous cete hypoth\`ese une autre module hilbertien $L^2_*$ et une
application holomorphe \'equivariante non constante vers un espace de Hilbert  $F: \wt{\X^{u}} \to L^2_*$. La $(1,1)$ forme $i\partial F \wedge\bar \partial F$ descend
au champ $\X$.  La classe de cohomolologie r\'esultante dans $H^{1,1}(\X, \R)=H^{1,1}(X,\R)$ est k\"ahl\'erienne par hypoth\`ese. 
Ceci observ\'e, la preuve donn\'ee par
\cite{Eys}, dans le cas non champ\^etre donne que $H\backslash\wt{\X^{u}}$ est de Stein, puis que $\wt{\X^{u}}$ est  Stein puisque le rev\^etement universel d'une vari\'et\'e de Stein est Stein.

\end{prv}

\subsection{Obstructions \`a la convexit\'e holomorphe} \label{obstr}

 Une obstruction \`a ce que $(\X, H)$ soit (HC) est la pr\'esence de chaines de Nori au sens suivant:

\begin{defi}\label{nori} Un sous champ complexe analytique ferm\'e de $\X$ est
 une chaine de Nori pour $(\X, H)$ si et seulement si les composantes connexes de
son image r\'eciproque dans $H\backslash \widetilde{\X^u}$ est une r\'eunion de composantes
irr\'eductibles compactes. 
\end{defi}

L'auteur ne connait pas d'exemples o\`u cette obstruction est effectivement pr\'esente. 
En revanche, un autre type d'obstruction est bien pr\'esent. Soit $A$ un tore complexe compact 
de dimension $g$ et
$\lambda:\pi_1(A)=\Z^{2g}\to \Z^r$ un morphisme surjectif.
On consid\`ere $(\X, H)=(A, \ker(\lambda))$.   Le th\'eor\`eme de
Remmert-Morimoto permet de d\'ecomposer $\ker(\lambda)\backslash \widetilde{A^u}$
comme un produit de la forme $\C^a\times\C^{*b}\times A'$ o\`u $A'$ est un
quotient de $\C^m$ n'ayant pas de fonction holomorphe non constante. Si $A'$ est
un tore complexe compact $\ker(\lambda)\backslash \widetilde{A^u}$ est holomorphiquement convexe,
mais si $A'$ est  non compact ce n'est trivialement pas le cas. Si $r<g$
et $A$ est simple (c'est \`a dire n'a pas de sous tore non compact), un tel quasi tore de Cousin est
automatiquement pr\'esent et la propri\'et\'e (HC) est viol\'ee. 

\begin{defi}
 Un compact lamin\'e $K\subset \X$ est de Cousin pour $(\X, H)$ si et seulement si $K$ est minimal et
il existe un voisinage ouvert $K\subset K^+$ tel que $\mathrm{Im}(\pi_1(K^+)\to \pi_1(\X))$
est fini. 
\end{defi}

\begin{lem}
 Si $(\X, H)$ est (HC) tout compact lamin\'e de Cousin est contenu dans un sous champ analytique
propre $Z\subset \X$ tel que $\mathrm{Im}(\pi_1(Z)\to \pi_1(\X))$
est fini. 
\end{lem}
\begin{prv}
La pr\'eimage de $K^+$ dans $H\backslash{\X^u}$ est une r\'eunion d'ouverts relativement compacts, 
donc il en va ainsi de la pr\'eimage de $K$.
Par minimalit\'e toute feuille de $K$ est dense et donc les fonctions holomorphes de  $H\backslash{\X^u}$ sont
constantes sur toute composante connexe $K'$ de la pr\'eimage de $K$. Par (HC) l'intersection des ensembles de niveaux
des fonctions holomorphes de $H\backslash{\X^u}$
d\'efinis  par leurs valeurs
sur $K'$ est un sous espace analytique compact de $H\backslash{\X^u}$ contenant $K'$. Son image
dans $\X$ est le sous champ
annonc\'e. 
\end{prv}

La conjecture de Shafarevich serait donc r\'efut\'ee si on pouvait exhiber  une vari\'et\'e projective-alg\'ebrique
de groupe fondamental  infini
avec   une chaine de Nori ou un compacts lamin\'e de Cousin
Zariski-dense.

\section{Le cas lin\'eaire}

La conjecture de Shafarevich est d\'emontr\'ee pour une vari\'et\'e K\"ahl\'erienne dont le groupe
fondamental est lin\'eaire, c'est \`a dire qu'il admet une repr\'esentation fid\`ele
dans $GL_N(\C)$, \cite{EKPR,CCE}. Nous donnons dans cette section un \'enonc\'e 
plus g\'en\'eral \'etendu au cadre champ\^etre. 

\subsection{Un r\'esultat g\'en\'eral de convexit\'e holomorphe}

Soit $G$ un groupe alg\'ebrique lin\'eaire r\'eductif sur $\Q$. 
Le groupe $\pi_1(\X)$ est de pr\'esentation finie.  On d\'efinit $R=R_B(\X, G)$ comme le sch\'ema des repr\'esentations
de  $\pi_1(\X)$ dans $G$, c'est \`a dire que $R$ repr\'esente le foncteur
 $$A\mapsto \mathrm{Hom_{Grp}}(\pi_1(\X), G(A))$$
$A$ parcourant la cat\'egorie des $\Q$-alg\`ebres. Il est connu que $R$ est affine de type fini. On note $\bar A [R]$ son anneau des fonctions r\'eguli\`eres
et on d\'esigne par $\rho_G^{\infty}: \pi_1(\X) \to G(\bar A [R])$ sa repr\'esentation tautologique. 

On note $M=M_B(\X, G):= R_B(\X, G) //G$ o\`u $G$ agit par conjugaison sur $R$. Le foncteur repr\'esent\'e par $M$ n'a pas la m\^eme expression simple que
pour $R$ mais les points de $M$ \`a valeurs dans un corps alg\'ebriquement clos de caract\'eristique nulle $\Omega$ sont les classes de conjugaisons
de l'action de $G(\Omega)$ sur $R(\Omega)^{ss}$ par conjugaison. La notation $R(\Omega)^{ss}$ d\'esigne l'ensemble des repr\'esentations 
\`a valeurs dans $G(\Omega)$ qui sont semisimples, i.e.: l'adh\'erence de Zariski de leur image est r\'eductive.
Pour chaque  composante irr\'eductible  $i$ de $M$ on note $\bar \rho^{0,i}_G: \pi_1(\X) \to G(\Omega_i)$ un repr\'esentant 
du point g\'en\'erique de $i$. On note  $\rho^i_G$   une des repr\'esentations complexes obtenues en  utilisant l'axiome du choix pour plonger le corps $\Omega_i$ dans $\C$.

\begin{theo}\label{main} Soit $\X$ un orbifold K\"ahl\'erien compact d\'eveloppable. 
On consid\`ere les sous groupes suivants $H^{\infty}_G=H^{\infty}_G(\X):=\ker (\rho_G^{\infty})$ et pour $I$ un ensemble fini non vide de composantes irr\'eductibles de $M$
$H^{0,I}_G=H^{0,I}_G(\X):=\bigcap_{i\in I} \ker(\rho^i_G)$. Alors $(\X, H_G^{\infty})$ et $(\X, H_G^{0,I})$ v\'erifient (HC). 
\end{theo}
\begin{prv}
 Les r\'esultats positifs g\'en\'eraux de \cite{Eys} \cite{EKPR} \cite{CCE} 
 \'etablissent cette conclusion
 si $\X$ est une vari\'et\'e k\"ahl\'erienne compacte.

Soit $K$ un corps et $\rho:\pi_1(\X)\to G(K)$ une repr\'esentation lin\'eaire d'image $\Gamma$. Par un th\'eor\`eme de
Malcev, il existe un sous groupe normal d'indice fini $H\triangleleft \pi_1(\X)$ tel  que $\rho(H)$
est sans torsion.  Soit $e:\X'\to \X$ le rev\^etement \'etale fini correspondant \`a $H$. 
L'image de la repr\'esentation $e^*\rho:\pi_1(\X')\to G(K)$ ne contient pas d'\'el\'ement de torsion non trivial. 
Soit $\gamma:\X'\to X'$ le morphisme naturel de $\X'$ sur son espace de modules.
Comme l'image par $\rho$ des groupes d'inertie de $\X'$ est triviale il suit du th\'eor\`eme
de Zariski-Van Kampen de
 \cite{No1} que $e^*\rho=\gamma^* \rho'$ o\`u $\rho':\pi_1(X')\to G(K)$ est une repr\'esentation. 
On peut de plus consid\'erer une r\'esolution des singularit\'es $\mu:\bar X'\to X'$
et construire $\mu^*\rho': \pi_1(\bar X')\to G(K)$ une repr\'esentation de groupe k\"ahl\'erien. 

Supposons que $K$ soit un corps local, possiblement non archim\'edien, et que $\rho$ est 
r\'eductive au sens de \cite{Eys}. La repr\'esentation $\mu^*\rho'$
est alors r\'eductive et on peut alors construire une application harmonique
$\mu^*\rho'$-\'equivariante comme comme dans $\cite{Eys}$. 
Cette application contracte les fibres de $\mu$ et donc est d\'efinie sur le rev\^etement universel 
de $X'$ donc sur celui de $\X'$. Par suite, on peut construire 
un fibr\'e de Higgs $\rho$-\'equivariant harmonique sur le rev\^etement universel de $\X$
dont la repr\'esentation d'holonomie est $\rho$
obtenant ainsi une g\'en\'eralisation orbifold de \cite{Sim2} dans le cas o\`u
$K$ est archim\'edien, de \cite{GrS} dans le cas non archim\'edien. On peut utiliser
 la d\'eformation $ (E,t\theta)_{t\to 0}$ du fibr\'e de Higgs 
$(E,\theta)$  \cite{Sim1,Sim2} pour conclure que toute repr\'esentation lin\'eaire de $\pi_1(\X)$ 
se d\'eforme \`a une variation de structure de Hodge complexe
  en raisonnant comme dans l'extension de ce r\'esultat au cas k\"ahl\'erien 
donn\'ee dans  \cite{CCE}. 

En utilisant \cite[Lemme 3.13, Lemme 5.3]{CCE} avec  $X:=\bar X'$  nous voyons que si $C$ est une composante
irr\'eductible de $M_B(\pi_1(\X), G)$,  les composantes irr\'eductibles de
sa pr\'eimage dans  $M_B(\pi_1(\X'), G)$ correspondent  \`a des constructibles absolus de $M_B(X', G)$. 
Ce point \'etant acquis, les arguments de \cite{CCE} permettent de conclure.

\end{prv}

On ne sait en revanche pas si pour toute repr\'esentation semisimple de $\pi_1(\X)$
le rev\^etement $\ker(\rho) \backslash \widetilde{\X^u}$ est (HC). Par \cite{CCE}, il n'y a pas de chaine 
de Nori au sens de la d\'efinition \ref{nori} pour $(\X, \ker(\rho))$. En revanche, des obstructions de Cousin
sont bel et bien pr\'esentes comme nous l'avons vu \`a la section \ref{obstr} ce qui n\'ecessite de prendre
des repr\'esentations lin\'eaires bien choisies ayant des noyaux  suffisamment gros.

\subsection{Morphisme de Shafarevich lin\'eaire}

Soit $G$ un groupe alg\'ebrique r\'eductifs d\'efini sur $\Q$ et $\mathfrak{s}_G: \X \to \mathfrak{S}_G(X)$
le morphisme de champs analytiques obtenu en appliquant le lemme \ref{shafmor}
 \`a $H=\ker(\rho^{\infty}_G)$ au vu du th\'eor\`eme
\ref{main}. En passant \`a l'espace de modules nous obtenons une fibration d'espaces analytiques normaux
$sh_G: X\to sh_G(\X)$ dont les fibres sont les sous espaces analytiques connexes $Z$ maximaux de $X$ tels que
$\rho_G^{\infty}(\pi_1(Z\times_X\X))$ est fini.

Le cas $G=GL_1$ est classique. En effet $sh_{GL_1}$ est exactement
 la factorisation de Stein de l'application d'Albanese de $\X$. 

Si $i:G\to G'$ est un morphisme de noyau fini
de groupes alg\'ebriques r\'eductifs les fibres de $sh_G$ contiennent celles de $sh_G'$ et on a 
une factorisation $sh_G'=sh_G \circ \sigma_{i}$. En particulier, la suite de fibrations holomorphes
$sh_{GL_N}: X\to sh_{GL_{N+1}}(\X)\to sh_{GL_N}(\X)$  stationne pour $N$ assez grand car
 $R_N:=X\times_{sh_{GL_N}(\X)}X$ forme une suite d\'ecroissante de ferm\'es analytiques propres de $X$.
On appelle morphisme de Shafarevich lin\'eaire la fibration
 $sh_{lin}: \X \to sh_{lin}(X)$ d\'efinie par $sh_{GL_N}$ pour $N\in \N$ assez grand.
Par construction, on a:
\begin{prop}
Les fibres de $sh_{lin}$ sont les les sous espaces analytiques connexes $Z$ maximaux de $X$ tels que
$\rho^{\infty}_{GL_N}(\pi_1(Z\times_X\X))$ est fini pour tout $N\in \N$.
\end{prop}

La dimension de Shafarevich lin\'eaire de $\X$ est par d\'efinition 
la dimension de  $sh_{lin}(X)$. Elle  est plus grande que la dimension d'Albanese de $\X$ et
s'y substitue avantageusement dans plusieurs 
applications. Par exemple, nous obtenons ais\'ement une 
am\'elioration du  r\'esultat classique de Napier \cite{Nap}:

\begin{theo}
 Soit $\X$ un orbifold k\"ahl\'erien compact d\'eveloppable de dimension $2$. On suppose que sa dimension
de Shafarevich lin\'eaire est $2$ et que $\X$ n'a pas de chaine de Nori au sens de la d\'efinition \ref{nori}. Alors,
 $\widetilde{\X^{u}}$ est holomorphiquement convexe.
\end{theo}
\begin{prv} La preuve de \cite{Nap} s'applique sans modification substantielle. 
 On notera qu'une \'eventuelle chaine de Nori serait support\'ee par une courbe exceptionnelle du morphisme $sh_{lin}$.
\end{prv}

Dans tous les exemples que l'auteur a pu analyser compl\`etement, il se trouve que 
l'application $\X\to sh_{lin}(\X)$  est un morphisme de Shafarevich pour $\X$, c'est \`a dire que
les fibres de $sh_{lin}$ sont les les sous espaces analytiques connexes $Z$ maximaux de $X$ tels que
$\mathrm{Im}(\pi_1(Z\times_X\X)\to \pi_1(\X))$ est fini. Lorsque cette propri\'et\'e plus forte
est v\'erifi\'ee, il d\'ecoule du th\'eor\`eme \ref{main} que $\widetilde{\X^{u}}$ est holomorphiquement convexe.
Cependant, \`a la connaissance de l'auteur, personne n'a pu exclure qu'un groupe k\"ahl\'erien
ait un compl\'et\'e profini infini (c'est \`a dire que l'intersection des noyaux des morphismes 
vers des groupes finis est un groupe d'indice infini) et \'egal \`a son compl\'et\'e proalg\'ebrique. 
Ceci se produirait si et seulement si toutes ses repr\'esentations lin\'eaires de dimension finie
\'etaient d'image finie.  Une telle vari\'et\'e serait de dimension de Shafarevich lin\'eaire nulle. 
Les exemples  de Toledo \cite{Tol} de vari\'et\'es projectives de groupe
fondamental non r\'esiduellement fini v\'erifient que $Sh_{lin}(X)=X$ ce qui 
implique que leur rev\^etement universel est de Stein. 

Finissons par une remarque. Puisque le morphisme d'Albanese d'une vari\'et\'e projective
 lisse poss\`ede une construction
purement alg\'ebrique, le morphisme $sh_{GL_1}$ \'egalement et on peut le d\'efinir sur un corps de base de
caract\'eristique positive. Nous ignorons compl\'etement si c'est le cas pour les $sh_{G}$ 
dont la construction est  transcendante. Une construction et une interpr\'etation
purement alg\'ebrique serait sans doute tr\`es int\'eressante, par exemple pour des questions d'arithm\'etique.
On peut notamment
esp\'erer qu'une vari\'et\'e alg\'ebrique $X$   d\'efinie sur un corps de nombres $K$ telle qu'il existe
$G$ avec $sh_G=\mathrm {id}_X$ v\'erifie $X(K)<+\infty$ si elle  ne contient pas de sous vari\'et\'e ab\'elienne.
Le cas $G=GL_1$ est un th\'eor\`eme classique de Faltings. 

\section{Observations relatives au cas singulier}

Une courbe rationnelle nodale a pour rev\^etement universel une
chaine infinie de composantes irr\'eductibles rationnelles et compactes. En particulier, il n'est pas
holomorphiquement convexe. Nous pensons que ceci est reli\'e au fait que l'homologie d'une telle courbe est pure de poids z\'ero.

Il semble cependant possible de trouver des familles $F$ de repr\'esentations du groupe
fondamental d'une vari\'et\'e k\"ahl\'erienne singuli\`ere $Z$ telle que 
$$\widetilde{Z_F}=\cap_{\rho \in F} \ker(\rho) \backslash \widetilde{Z^u}$$ soit holomorphiquement convexe.
Par exemple, s'il existe $i:Z\to \X$ un morphisme vers un orbifold lisse et $\rho^{\infty}_G$
est comme dans le th\'eor\`eme \ref{main}, il suit de ce m\^eme th\'eor\`eme que
$\ker(i^*\rho^{\infty}_G) \backslash \widetilde{Z^u}$ est holomorphiquement convexe.
Ceci sugg\`ere l'existence d'une th\'eorie de Hodge Mixte pour le $H^1$ non ab\'elien $H^1(Z,G)$
de sorte que les sous ensembles de repr\'esentations correspondant une sous structure de Hodge pure 
de poids 1
forment une famille $F$ de repr\'esentations de $\pi_1(Z)$ telle que 
$\cap_{\rho \in F} \ker(\rho) \backslash \widetilde{Z^u}$ soit holomorphiquement convexe.
Malheureusement, nous ne savons pas d\'evelopper les \'el\'ements n\'ecessaires d'une
 th\'eorie de Hodge mixte non ab\'elienne.

Voici toutefois un r\'esultat dans cette direction que nous avons utilis\'e sans donner de d\'emonstration dans \cite{CCE}. 

\begin{theo}
 Soit $Z$ un espace complexe k\"ahl\'erien compact  et $H\subset H^1(Z,\Q)$
une sous structure de Hodge pure de poids $1$ de la structure de Hodge mixte de Deligne
sur $H^1(Z,\Q)$. Alors il existe une application holomorphe de la seminormalisation  $sn:Z^{sn} \to Z$
 \footnote{Noter que $sn$ est un hom\'eomorphisme et que $sn_*O_{Z^{sn}}$ est le sous faisceau du faisceau des
fonctions continues de $Z$ qui sont holomorphes sur la partie lisse de $Z$.}
vers un tore complexe $a:Z^{sn}\to T$ telle que $H$ soit l'image de $a^*: H^1(T,\Q)\hookrightarrow 
 H^1(Z^{sn},\Q)=H^1(Z^{sn},\Q)$. 

En particulier,
$\cap_{\rho \in H} \ker(\rho) \backslash \widetilde{Z^u}$ est holomorphiquement convexe.
\end{theo}
\begin{prv}
 Soit $S$ un tore complexe compact. Consid\'erons deux points distincts $x, y\in S$ et construisons
une vari\'et\'e k\"ahl\'erienne seminormale $S_{x\sim y}$ en identifiant $x$ et $y$. On a un morphisme
naturel $S\to S_{x\sim y}$ qui donne lieu a une extension de structures de Hodge mixtes enti\`eres:
$$ 0\to \Z (0) \to H^1(S_{x\sim y}) \to H^1(S,\Z) \to 0.
$$
Le groupe des extensions de structures de Hodge mixtes enti\`eres de $H^1(S,\Z)$ par $\Z (0)$ est
exactement $S$ et l'extension pr\'ec\'edente a pour classe l'invariant d'Abel Jacobi $x-y\in S$. 
Par suite, cette extension n'est jamais scind\'ee. 

Plus g\'en\'eralement,  l'extension est scind\'ee sur $\Q$
si et seulement si $x-y$ est un point de torsion.

Soit $T$ un tore complexe tel que $H^1(T,\Q)=H$ comme $\Q$-structure de Hodge. 
Quitte a prendre une isog\'enie de $T$ on peut supposer l'existence d'un morphisme de SHM $H^1(T,\Z)\to 
H^1(Z, \Z)$.  

Supposons d'abord que $Z$ est normale. Dans ces conditions la structure de Hodge de $H^1(Z)$ est pure et le 
r\'esultat d\'ecoule de la th\'eorie du morphisme d'Albanese. 

Supposons maintenant que la normalisation $Z^{no}$ de $Z$ est irr\'eductible. Consid\'erons 
le morphisme $a:Z^{no}\to T$ construit au pas pr\'ec\'edent.  Fixons deux points distincts $p,q$ de $Z^{no}$
ayant m\^eme image  dans $Z$. On a un diagramme $Z^{no}\to Z^{no}_{p\sim q}\to Z^{sn}\to Z$. 

On dispose d'un diagramme dont les lignes horizontales sont exactes:
$$
\begin{array}{ccccccccc}
 0 & \to & \Z(0) & \to &H^1(T_{a(p)\sim a(q)},\Z) &\to& H^1 (T, \Z) & \to& 0\\
& & \downarrow  & & \downarrow & & \downarrow & & \\
 0 & \to & \Z(0) & \to &H^1(Z^{no}_{p\sim q},\Z) &\to& H^1 (Z^{no}, \Z) & \to& 0
\end{array}
$$

Ceci identifie comme $\Z$-structure de Hodge mixte \`a $H^1(T_{a(p)\sim a(q)},\Z)$ la pr\'eimage de $H^1(T, \Z)$
dans $H^1 ( Z^{no}_{p\sim q},\Z)$. Mais le morphisme de SHM $H^1(T, \Z)\to H^1(Z^{no}, \Z)$ factorise
comme $H^1(T,\Z)\to H^1(Z, \Z) \to H^1(Z^{no},\Z)$ et on a un morphisme de SHM 
$H^1(Z, \Z)\to H^1(Z^{no}_{p\sim q},\Z)$ qui permet de construire un morphisme de SHM
$H^1(T,\Q) \to H^1(Z^{no}_{p\sim q},\Z)$ scindant $H^1(T_{a(p)\sim a(q)},\Z)$. Donc $a(p)-a(q)$ est nul. 
Ceci implique que l'application $a$
descend ensemblistement \`a $Z$ et on voit ais\'ement qu'elle y est continue. Ceci donne bien une application
holomorphe $Z^{sn}\to T$. 

Dans le cas g\'en\'eral $Z=Z_1\cup Z_2$ avec $Z_1$ et $Z_2$ non vides connexes. Si $Z_1\cap Z_2=\{o\}$ est 
compos\'e d'un seul point on recolle $a|_{Z_1}$ avec $a|_{Z_2}$ en translatant la seconde de mani\`ere \`a ce que 
l'application soit uniqument d\'efinie en $o$. Cela est coh\'erent car
 $H^1(Z)=H^1(Z_1)\oplus H^1(Z_2)$ dans le cas pr\'esent. 
Dans le cas g\'en\'eral, on choisit deux points base de $Z_1$ et $Z_2$ s'identifiant dans $Z$, on 
construit $Z'$ en identifiant ces deux points bases et le raisonnement pr\'ec\'edent s'applique 
pour voir que deux points de $Z'$ s'identifiant dans $Z$ ont m\^eme image par $a:Z'\to T$. 

\end{prv}

{Philippe Eyssidieux}\\
{Institut Universitaire de France et Universit\'e  Grenoble I. Grenoble, France}\\
{philippe.eyssidieux@ujf-grenoble.fr}\\
{http://www-fourier.ujf-grenoble.fr/$\sim$eyssi/}
\end{document}